\documentclass{article}
\usepackage{amsmath,amsfonts,amssymb}

\DeclareMathSymbol\nullset{\mathord}{AMSb}{"3F}

\begin{document}

\markboth{C. Bekh-Ochir and S. A. Rankin}
{Central polynomials}

\title{The central polynomials of the infinite dimensional unitary and nonunitary Grassmann algebras}

\author{C. Bekh-Ochir and S. A. Rankin}

\maketitle

\begin{abstract}
We describe the $T$-space of central polynomials for both the unitary 
and the nonunitary infinite dimensional Grassmann algebra over a field 
of characteristic $p\ne2$ (infinite field in the case of the unitary algebra).
\end{abstract}



\newcounter{parts}
\def\set#1\endset{\{\,#1\,\}}
\def\gen#1{\left<{#1}\right>}
\def\rest#1{{}_{{}_{#1}}}
\def\com#1,#2{[{#1},{#2}]}
\def\choice#1,#2{\binom{#1}{#2}}
\def\kx{k\langle X\rangle}
\def\kzerox{k_0\langle X\rangle}
\def\konex{k_1\langle X\rangle}
\def\unitgrass{G}
\def\nonunitgrass{G_0}
\def\finiteunitgrass#1{G(#1)}
\def\finitenonunitgrass#1{G_0(#1)}
\def\siderovset{SS}
\def\finitesiderovsetone#1{SS'(#1)}
\def\finitesiderovsettwo#1{SS''(#1)}
\def\boundedsiderovset{BSS}
\def\bss#1{BSS(#1)}
\def\lend#1{lend(#1)}
\def\lbeg#1{lbeg(#1)}
\def\finitesidset#1{SS(#1)}
\let\cong=\equiv
\def\mod#1{\,\,(\text{mod}\,#1)}

\newtheorem{theorem}{Theorem}[section]
\newtheorem{corollary}{Corollary}[section]
\newtheorem{lemma}{Lemma}[section]
\newtheorem{proposition}{Proposition}[section]
{\newtheorem{definition}{Definition}[section]}
\def\proof{\ifdim\lastskip<\smallskipamount\relax\removelastskip
  \vskip\smallskipamount\fi\leavevmode\noindent\hbox to 0pt{\hfil}{\it Proof.}}
\def\strutdepth{\dp\strutbox}
\def\epmarker{\vbox to \strutdepth{\baselineskip\strutdepth\vss\hfill{%
\hbox to 0pt{\hss\vrule height 4pt width 4pt depth 0pt}\null}}}
\def\edproofmarker{\strut\vadjust{\kern-2\strutdepth\epmarker}}
\def\endproof{\edproofmarker\vskip10pt}


\section{Introduction and preliminaries}
 Let $k$ be a field and $X$ a countable set, say $X=\set x_i\mid i\ge 1\endset$. Then $\kzerox$ denotes
 the free (nonunitary) associative $k$-algebra over $X$, while $\konex$ denotes the free unitary associative
 $k$-algebra over $X$. 

 Let $H$ denote any associative $k$-algebra. 
 For any $X\subseteq H$, $\gen{X}$ shall denote the linear subspace of $H$ spanned by $X$.
 Any linear subspace of $H$ that is invariant under every endomorphism of $H$ is
 called a $T$-space of $H$, and if a $T$-space happens to also be an ideal of $H$, then it
 is called a $T$-ideal of $H$. For $X\subseteq H$, the smallest $T$-space containing $X$
 shall be denoted by $X^S$, while the smallest $T$-ideal of $H$ that contains $X$ shall be
 denoted by $X^T$. In this article, we shall deal only with $T$-spaces and
 $T$-ideals of $\kzerox$ and $\konex$.
 
 An element $f\in \kzerox$ is called an {\em identity} of $H$
 if $f$ is in the kernel of every homomorphism from $\kzerox$ to $H$ (from $\konex$ if $H$ is unitary). 
 The set of all identities of $H$ is a $T$-ideal of $\kzerox$ (and of $\konex$ if $H$ is unitary), 
 denoted by $T(H)$. An element $f\in \kzerox$ is called a {\it central polynomial} of $H$ 
 if $f\notin T(H)$ and the image of $f$ under any homomorphism from $\kzerox$ ($\konex$ if $H$ is unitary)
 belongs to $C_{H}$, the centre of $H$.
 
 Let $\unitgrass$ denote the (countably) infinite dimensional unitary
 Grassmann algebra over $k$, so there exist $e_i\in \unitgrass$, $i\ge 1$, 
 such that  for all $i$ and $j$, $e_ie_j=-e_je_i$, $e_i^2=0$, and
 $\mathcal{B}=\set e_{i_1}e_{i_2}\cdots e_{i_n}\mid n\ge 1,
 i_1<i_2<\cdots i_n\endset$, together with $1$, forms a linear basis for
 $G$. The subalgebra of $\unitgrass$ with linear basis $\mathcal{B}$ is
 the infinite dimensional nonunitary Grassmann algebra over $k$, and is
 denoted by $\nonunitgrass$.  

 Evidently, $T(\unitgrass)\subseteq T(\nonunitgrass)$. It is well known
 that $T^{(3)}$, the $T$-ideal of $\konex$ that is generated  by $\com
 {\com x_1,{x_2}},{x_3}$, is contained in $T(\unitgrass)$. For convenience, 
 we shall write $\com x_1,{x_2,x_3}$ for  $\com {\com x_1,{x_2}},{x_3}$.
 Moreover, it is also well known (see \cite{D}, Theorem 5.1.2 (i) for example)
 that for a field of characteristic zero, $T(G)=T^{(3)}$.

 We shall let $CP(\unitgrass)$ and $CP(\nonunitgrass)$ denote
 the $T$-spaces of $\konex$ and $\kzerox$ that are generated by the 
 central polynomials of $\unitgrass$ and $\nonunitgrass$, respectively. 
 Evidently, $T(\unitgrass)\subseteq CP(\unitgrass)$, $T(\nonunitgrass)
 \subseteq CP(\nonunitgrass)$, and $CP(\unitgrass)\cap \kzerox\subseteq CP(\nonunitgrass)$.

 The original interest in central polynomials dates back to 1956 \cite{Ka}, when Kaplansky 
 asked whether there exists a central polynomial for the matrix algebra $M_r(k)$, $r>2$.
 It wasn't until 1972-73 that the answer to this question was provided by
 constructions due to Formanek \cite{Fo} and Razmyslov \cite{Raz}.
 Subsequently, many authors have contributed examples of central polynomials for
 the matrix rings.

 In contrast, little seems to be known about the central polynomials for another natural
 algebra to study, namely the Grassmann algebra.
 It is known (see \cite{Br}) that for a field of characteristic zero, $CP(G)$, the $T$-space
 of central polynomials of $G$, the infinite dimensional unitary Grassmann algebra,
 is equal to $\set \com x_1,{x_2}, \com x_1,{x_2}\com x_3,{x_4}\endset^S$.  Furthermore, for
 a field of characteristic 2, the unitary infinite dimensional Grassmann
 algebra (and hence the nonunitary Grassmann algebra) is commutative,
 and thus the $T$-space of central polynomials for each is $\kzerox$.

 In this paper, we describe the $T$-space of the central polynomials for
 both the unitary and the nonunitary infinite dimensional Grassmann
 algebra over a field of characteristic $p\ne 2$ (except for the case of
 the infinite dimensional unitary Grassmann algebra over a finite field).
 As should be expected, the description of the central polynomials relies heavily on
 knowledge of the $T$-ideal of identities of the
 relevant Grassmann algebra, due to A. Giambruno and P.
 Koshlukov \cite{Gi} for the case of the unitary infinite dimensional
 Grassmann algebra over an infinite field, and to P. N. Siderov
 \cite{Si} for the case of the nonunitary infinite dimensional Grassman
 algebra over an arbitrary field of characteristic $p>2$. 

 In the case of the infinite dimensional nonunitary Grassmann algebra over an
 arbitrary field $k$, we show that
 if $k$ has characteristic zero, then 
 $$
  CP(\nonunitgrass)=\set \com x_1,{x_2},\com x_1,{x_2}\com x_3,{x_4}\endset^S,
 $$
 while if $k$ has characteristic $p>2$, then
 $$
   CP(\nonunitgrass)=\set \com x_1,{x_2},x_1^p,x_2x_1^p,\com x_1,{x_2}\com x_3,{x_4}\endset^S+\set w_n\mid  n\ge 1\endset^S,
 $$   
 where for each $n\ge1$, $w_n=\prod_{k=1}^n\com {x_{2k-1}},{x_{2k}}x_{2k-1}^{p-1}x_{2k}^{p-1}$.

 In the case of the infinite dimensional unitary Grassman algebra over an infinite
 field $k$, we show that if $k$ has characteristic zero, then 
 $$
  CP(\unitgrass)=
 \set \com x_1,{x_2}, \com x_1,{x_2}\com x_3,{x_4}\endset^S,
 $$
 while if $k$ has characteristic $p>2$, then
 $$
  CP(\unitgrass)=\set \com x_1,{x_2},x_1^p,\com x_1,{x_2}\com x_3,{x_4}\endset^S+\set x_{2n+1}^pw_n\mid  n\ge 1\endset^S,
 $$
 where for each $n\ge1$, $w_n=\prod_{k=1}^n\com {x_{2k-1}},{x_{2k}}x_{2k-1}^{p-1}x_{2k}^{p-1}$.
 
 We complete this section with a description of some familiar identities 
 that will be required in the sequel.

\begin{lemma}\label{lemma: handy}
 \begin{list}{(\roman{parts})}{\usecounter{parts}}
  \item 
    $\com u,{vw}=\com u,vw+v\com u,w$ for all $u,v,w\in \kzerox$.
  \item 
    $\com u,{vw}=\com u,vw+\com u,wv+\com v,{\com u,w}$ for all $u,v,w\in \kzerox$.
  \item 
     $\com u,{\prod_{i=1}^n v_i}\cong \sum_{i=1}^n 
          \bigl(\mkern3mu \com u,{v_i}\prod_{\substack{j=1\\j\ne i}}^n v_j\mkern3mu\bigr)\mod{T^{(3)}}$ for any positive 
	  integer $n$, and any $u,v_1,v_2,\ldots,v_n\in \kzerox$. 
   \item  
     $\com u,v\com w,x\cong -\com u,w\com v,x\mod{T^{(3)}}$ for all $u,v,w\in \kzerox$.
   \item  
    $\com u,v\com u,w\cong 0\mod{T^{(3)}}$ for all $u,v,w\in \kzerox$.
   \item  
    $\com u,v uw\cong \com u,vwu\mod{T^{(3)}}$ for all $u,v,w\in \kzerox$.    
  \item 
   For $n\ge2$, $x_1^nx_2^n\cong(x_1x_2)^n+\choice n,2\com x_1,{x_2}x_1^{n-1}x_2^{n-1}\mod{T^{(3)}}$.    
 \end{list}
\end{lemma}

\begin{proof}
 Parts (i) and (ii) are evident by direct calculation. We prove (iii) by induction on $n$. The result is trivial 
 when $n=1$, while  (ii) establishes the case for $n=2$.
 Suppose now that $n\ge2$ is an integer for which the result holds, and consider $u,v_1,v_2,\ldots,v_n,v_{n+1}\in\kzerox$.
 Working modulo $T^{(3)}$, we have
 \begin{align*}
  \com u,{\prod_{i=1}^{n+1} v_i}&\cong\com u,{\prod_{i=1}^n v_i} v_{n+1} + \com u,{v_{n+1}}\prod_{j=1}^n v_i\quad\mbox{by (i)}\\
                                &\cong\sum_{i=1}^n\bigl(\mkern3mu\com u,{v_i}\prod_{\substack{j=1\\j\ne i}}^{n+1}v_i\mkern3mu\bigr)+\com u,{v_{n+1}}\prod_{j=1}^n v_i\quad\text{by the inductive hypothesis}\\
				&=\sum_{i=1}^{n+1}\bigl(\mkern3mu\com u,{v_i}\prod_{\substack{j=1\\j\ne i}}^{n+1}v_i\mkern3mu\bigr).
  \end{align*}

 The result follows now by induction.

 For (iv) and (v), see Latyshev \cite{L}. 
 (vi) follows immediately from (v), as we have $\com u,v uw-\com u,v wu=\com u,v\com u,w\cong 0\mod{T^{(3)}}$.
 
 For (vii), observe that $x_1^2x_2^2=(x_1x_2)^2+x_1\com x_1,{x_2}x_2\cong (x_1x_2)^2+\com x_1,{x_2}x_1x_2$
 modulo $T^{(3)}$. Suppose now that $n\ge 2$ is an integer such that $x_1^nx_2^n\cong(x_1x_2)^n+\choice n,2\com x_1,{x_2}x_1^{n-1}x_2^{n-1}$,
 again modulo $T^{(3)}$.
 Then 
 \begin{align*}
    x_1^{n+1}x_2^{n+1}&=x_1(x_1^nx_2^n)x_2\cong x_1(x_1x_2)^nx_2+\choice n,2x_1\com x_1,{x_2}x_1^{n-1}x_2^{n}\\
     &\cong (x_1x_2)^{n+1}+\com x_1,{(x_1x_2)^n}x_2+\choice n,2\com x_1,{x_2}x_1^nx_2^n\\
     &\cong (x_1x_2)^{n+1}+n\com x_1,{x_1x_2}(x_1x_2)^{n-1}x_2+\choice n,2\com x_1,{x_2}x_1^nx_2^n\quad\text{by (iii)}\\
     &\cong (x_1x_2)^{n+1}+\choice n+1,2\com x_1,{x_2}x_1^nx_2^n\quad\text{by (vi)},
 \end{align*}
 and so the result follows by induction.
\end{proof}

We note that by definition, any 2-commutator is central modulo $T^{(3)}$, while Lemma \ref{lemma: handy} (vi)
establishes that in any product term that has a
2-commutator factor $\com u,v$, any occurrence of $u$ (or $w$) in the product term commutes with any factor of the product term, 
modulo $T^{(3)}$.


\section{The central polynomials of the infinite dimensional nonunitary Grassmann algebra over a field of characteristic $p\ne 2$}
Throughout this section, it is assumed that $k$ has characteristic $p\ne2$.

\begin{lemma}\label{lemma: 2 commutator commutes}
 For any $u,v\in \kzerox$ and nonnegative integers $m,n$, $\com u,v u^mv^n$ 
 is central modulo $T^{(3)}$.
\end{lemma}

\begin{proof}
 First of all, for any $w\in\konex$, $\com w,{\com u,v}=w\com u,v - \com u,vw$, so $\com u,v$ is indeed
 central modulo $T^{(3)}$. With this fact and Lemma \ref{lemma: handy} (vi), we then obtain
 $w\com u,v u^mv^n\cong \com u,v u^mv^nw\mod{T^{(3)}}$.
\end{proof}
 
\begin{definition}\label{definition: def of S}
 If $p=0$, let $S=\set\com x_1,{x_2}\endset^S$, while if $p>2$, then for each $n\ge1$, let 
 $w_n=\prod_{k=1}^n\com {x_{2k-1}},{x_{2k}}x_{2k-1}^{p-1}x_{2k}^{p-1}$ and set
 $S=\set \com x_1,{x_2}\endset^S+\set w_n\mid \ n\ge 1\endset^S$.
\end{definition}

 \begin{lemma}\label{lemma: s central}
  $S\subseteq CP(\unitgrass)\cap\kzerox\subseteq CP(\nonunitgrass)$.
 \end{lemma}
 
 \begin{proof}
  This is immediate from Lemma \ref{lemma: 2 commutator commutes}.
 \end{proof} 

\begin{lemma}\label{lemma: need for alt description of S gen}
 Let $\alpha_1,\alpha_2$ be positive integers, and suppose that $u\in\kzerox$ is central modulo $T^{(3)}$.
 Then $\com x_2,{x_1^{\alpha_1+1}x_2^{\alpha_2}u}\cong (\alpha_1+1)\com x_2,{x_1} x_1^{\alpha_1}x_2^{\alpha_2}u\mod{T^{(3)}}$.
\end{lemma}

\begin{proof}
 Working modulo $T^{(3)}$, we have
 \begin{align*}
  \com x_2,{x_1^{\alpha_1+1}x_2^{\alpha_2}u}&\cong \com x_2,{x_1^{\alpha_1+1}x_2^{\alpha_2}}u+
              \com x_2,u x_1^{\alpha_1+1}x_2^{\alpha_2} \quad{\text{by Lemma \ref{lemma: handy} (ii)}}\\ 
   &\cong  \com x_2,{x_1^{\alpha_1+1}x_2^{\alpha_2}}u
               \quad{\text{since $\com x_2,u\cong0\mod{T^{(3)}}$}}\\ 
   &\cong  \com x_2,{x_2^{\alpha_2}}x_1^{\alpha_1+1}u+\com x_2,{x_1^{\alpha_1+1}}x_2^{\alpha_2}u \quad{\text{by Lemma \ref{lemma: handy} (ii)}}\\ 
   &\cong (\alpha_1+1)\com x_2,{x_1} x_1^{\alpha_1}x_2^{\alpha_2}u\quad{\text{by Lemma \ref{lemma: handy} (iii).}}
 \end{align*}
\end{proof}
  
\begin{corollary}\label{corollary: alt description of S gen}
 Let $\beta_i$, $i=1\ldots,2n$ be integers such that for each $i$, $0\le \beta_i$, and, if $p>2$, $\beta_i\le p-1$. Then 
 \[
  \prod_{k=1}^n\com {x_{2k-1}},{x_{2k}}x_{2k-1}^{\beta_{2k-1}}x_{2k}^{\beta_{2k}}\in S+T^{(3)}.
 \] 
\end{corollary}

\begin{proof}
 Let $u=\prod_{k=1}^n\com {x_{2k-1}},{x_{2k}}x_{2k-1}^{\beta_{2k-1}}x_{2k}^{\beta_{2k}}$.
 If $p>2$ and $\beta_i=p-1$ for each $i$, then $u\in S$, while if either $p>2$ and there exists $i$ such that $\beta_i<p-1$, or
 $p=0$ and $i$ is arbitrary, then $\beta_i+1$ is invertible
 and so, working modulo $T^{(3)}$ (by Lemma \ref{lemma: handy} (vi), we may assume without loss of generality that $i$ is odd, say $i=2j-1$), we have
 \[
   u\cong-(\beta_i+1)^{-1}\com x_{i+1},{x_i^{\beta_i+1}x_{i+1}^{\beta_{i+1}}\bigl(\prod_{\substack{1\le k\le n\\ k\ne j}} \com {x_{2k-1}},{x_{2k}}x_{2k-1}^{\beta_{2k-1}}x_{2k}^{\beta_{2k}}\bigr)}\in S\\
 \]  
  by Lemma \ref{lemma: 2 commutator commutes} and Lemma \ref{lemma: need for alt description of S gen}.  
\end{proof}

\begin{corollary}\label{corollary: need for R gen}
 Let $\alpha_1,\alpha_2$ be nonnegative integers, let $u,v\in\kzerox$, and suppose that $u\in\kzerox$ is central modulo $T^{(3)}$.
 Then $(\alpha_1+1)v\com x_1,{x_2} x_1^{\alpha_1}x_2^{\alpha_2}u\cong\com x_2,v{x_1^{\alpha_1+1}x_2^{\alpha_2}u} \mod{S+T^{(3)}}$.
\end{corollary}

\begin{proof}
  By Lemma \ref{lemma: handy} (ii), we have
  \[
  \com x_2,{vx_1^{\alpha_1+1}x_2^{\alpha_2}u}\cong v\com x_2,{x_1^{\alpha_1+1}x_2^{\alpha_2}u}+
  \com x_2,v x_1^{\alpha_1+1}x_2^{\alpha_2}u\mod{T^{(3)}},
  \]
 so  $\com x_2,v x_1^{\alpha_1+1}x_2^{\alpha_2}u\cong -v\com x_2,{x_1^{\alpha_1+1}x_2^{\alpha_2}u}\mod{S+T^{(3)}}$,
 and by Lemma \ref{lemma: need for alt description of S gen}, 
 \[
  \com x_2,{x_1^{\alpha_1+1}x_2^{\alpha_2}u}\cong
 -(\alpha_1+1)\com x_1,{x_2} x_1^{\alpha_1}x_2^{\alpha_2}u\mod{T^{(3)}}.
 \] 
 Since $T^{(3)}$ is an ideal of $\kzerox$, the result follows.
\end{proof}
 
The following proposition is due to Siderov (presented in \cite{Si} for the nonunitary case).

\begin{proposition}\label{proposition: 2-comm ideal}
 Every element of $\set \com x_1,{x_2}\endset^T$ is congruent modulo $T^{(3)}$ to a linear combination of elements of the form
 $w\prod_{s=1}^n \com x_{j_{2s-1}},{x_{2s}}x_{j_{2s-1}}^{\beta_{2s-1}}x_{j_{2s}}^{\beta_{2s}}$, 
 where $j_1<j_2<\cdots j_{2s}$, $\beta_i\ge0$ for all $i$, and either $w$ is void, or else
 $w=\prod_{r=1}^m x_{i_r}^{\alpha_r}$ and $i_1<i_2<\cdots < i_r$, 
 $\set i_1,\ldots,i_r\endset\cap\set j_1,\ldots,j_{2s}\endset=\nullset$, and $\alpha_i>0$ for all $i$.
\end{proposition}

\begin{proof}
 Let $f\in\set \com x_1,{x_2}\endset^T$, so that $f$ is a linear combination of elements of the form $w\com u,v z$, where
 $w,u,v,z$ are monomials in $\kzerox$, with one or both of $w,z$ possibly void. 
 Now for monomials $u$ and $v=\prod_{s=1}^n x_{j_s}$, we may apply Lemma \ref{lemma: handy} (iii) to obtain that
 \[
   \com u,v\cong \sum_{s=1}^n \com u,{x_{j_s}}\prod_{\substack{t=1\\t\ne s}}^n x_{j_t}\mod {T^{(3)}}.
 \]
 Additionally, since $2$-commutators commute modulo $T^{(3)}$, it follows that $f$ can be
 written as a linear combination of elements of the form
 $w\prod_{s=1}^n \com x_{j_{2s-1}},{x_{j_{2s}}}$, 
 where either $w$ is void, or else
 $w$ has the form $w=\prod_{r=1}^m x_{i_r}^{\alpha_r}$ with $\alpha_l>0$ for each $l$. 
 Consider a polynomial of the form $w\com u,v$, where $w=\prod_{r=1}^m x_{i_r}^{\alpha_i}$ with 
 $\alpha_i>0$ for each $i$, and suppose that the variables appearing in $w$ are not in the proper order. 
 Then we use the identity $xy=yx+\com x,y$ to replace $w$ by a sum of two terms, one a monomial on the same variables as $w$ but
 with fewer out of order, the other of the form $w'[x_i,x_j][u,v]$, where $w'$ has fewer variables than $w$. It follows
 that $x$ can be written as a linear combination of polynomials of the
 form $w\prod_{s=1}^n \com x_{j_{2s-1}},{x_{2s}}$, where 
 $w$ is either void or else $w=\prod_{r=1}^m x_{i_r}^{\alpha_r}$, where $i_1<i_2<\cdots < i_r$, and 
 $\alpha_i> 0$ for all $i$. Next, we use Lemma \ref{lemma: handy} (iv) and (vi) to obtain $x$ as a linear combination of
 elements of the form
 $w\prod_{s=1}^n \com x_{j_{2s-1}},{x_{2s}}x_{j_{2s-1}}^{\beta_{2s-1}}x_{j_{2s}}^{\beta_{2s}}$, 
 where $j_1<j_2<\cdots j_{2s}$, $\beta_l\ge0$ for all $l$, and either $w$ is void, or else
 $w=\prod_{r=1}^m x_{i_r}^{\alpha_r}$ and $i_1<i_2<\cdots < i_r$, 
 and $\alpha_i>0$ for all $i$. Finally, we use Lemma \ref{lemma: handy} (vi) to obtain $f$ as a linear combination of
 elements of the desired form. 
\end{proof}

The elements that appear in Proposition \ref{proposition: 2-comm ideal} will play an important role in the sequel.

\begin{definition}\label{definition: siderov's elements}
 Let $\siderovset$ denote the set of all elements of the form 
 \begin{list}{(\roman{parts})}{\usecounter{parts}}
  \item $\prod_{r=1}^t x_{i_r}^{\alpha_r}$, or
  \item $\prod_{r=1}^s \com x_{j_{2r-1}},{x_{2r}}x_{j_{2r-1}}^{\beta_{2r-1}}x_{j_{2r}}^{\beta_{2r}}$, or
  \item $\bigl(\prod_{r=1}^t x_{i_r}^{\alpha_r}\bigr)\prod_{r=1}^s \com x_{j_{2r-1}},{x_{2r}}x_{j_{2r-1}}^{\beta_{2r-1}}x_{j_{2r}}^{\beta_{2r}}$, 
 \end{list}  
 \noindent where $j_1<j_2<\cdots j_{2s}$, $\beta_i\ge0$ and $\alpha_i\ge 1$ for all $i$, $\set i_1,\ldots,i_r\endset
 \cap\set j_1,\ldots,j_{2s}\endset=\nullset$, and $i_1<i_2<\cdots < i_t$.

 Let $u\in \siderovset$. If $u$ is of the form (i), then the beginning of $u$ is
 $\prod_{r=1}^t x_{i_r}^{\alpha_r}$, the end of $u$ is empty, the length of the beginning of $u$, $\lbeg{u}$, is
 equal to $t$ and the length of the end of $u$, $\lend{u}$, is 0. If $u$ is of the form
 (ii), then we say that the beginning of $u$ is empty, the end of $u$ is $\prod_{r=1}^s \com x_{j_{2r-1}},{x_{2r}}x_{j_{2r-1}}^{\beta_{2r-1}}x_{j_{2r}}^{\beta_{2r}}$,
 and $\lbeg{u}=0$ and $\lend{u}=s$.
 If $u$ is of the form (iii), 
 then we say that the beginning of $u$ is $\prod_{r=1}^t x_{i_r}^{\alpha_r}$, the end of $u$ is
 $\prod_{r=1}^s \com x_{j_{2r-1}},{x_{2r}}x_{j_{2r-1}}^{\beta_{2r-1}}x_{j_{2r}}^{\beta_{2r}}$, $\lbeg{u}=t$ and
 $\lend{u}=s$.
\end{definition}

In \cite{At}, Venkova introduced a total order on the set $\siderovset$ which was useful in her work on the
identities of the finite dimensional nonunitary Grassmann algebra. Subsequently, Siderov \cite{Si} used a modification
of this order in his work on the identities of the infinite dimensional nonunitary Grassmann algebra.

\begin{definition}[Venkova's ordering]\label{definition: total order}
 For $u,v\in \siderovset$, we say that $u>v$ if one of the following requirements holds.
  \begin{list}{(\roman{parts})}{\usecounter{parts}}
  \item $\deg u < \deg v$.
  \item $\deg u = \deg v$ but $\lend{u}<\lend{v}$.
  \item $\deg u = \deg v$ and $\lend{u}=\lend{v}$, but there exists $i \ge 1$ such that 
        $\deg_{x_i} u < \deg_{x_i} v$ and for each $j<i$, $\deg_{x_j} u = \deg_{x_j} v$.
  \item $\deg u = \deg v$, $\lend{u}=\lend{v}$ and for each $i \ge 1$,
        $\deg_{x_i} u = \deg_{x_i} v$, and there exists $j \ge 1$ such that 
        $x_j$ appears in the end of $u$ and in the beginning of $v$, and 
        for each $k<j$, $x_k$ appears in the beginning of $u$ if and only if $x_k$ appears in 
        the beginning of $v$.
  \end{list}
 \end{definition}
 
 It will be helpful to note that if $u>v$ by virtue of condition (iv), then there exists $k>j$ such that $x_k$ is in
 the beginning of $u$ and in the end of $v$.

\begin{definition}
 Let $R$ be the subspace of $\kzerox$ that is spanned by $\set u\in \siderovset\mid \lbeg{u}>0\endset$.
\end{definition} 

\begin{lemma}\label{lemma: rep of kzerox}
 $\kzerox=R+S+T^{(3)}$.
\end{lemma}

\begin{proof}
 Observe that by Proposition \ref{proposition: 2-comm ideal} and Corollary \ref{corollary: alt description of S gen},
 each element of $\set \com x_1,{x_2}\endset^T$ 
 is in $R+S+T^{(3)}$, and so it suffices to prove that every monomial belongs to $R+S+T^{(3)}$. This
 we prove by induction on the number of pairs of indices that are out of order.
 Consider an arbitrary monomial $w=x_{t_1}^{\alpha_1}x_{t_2}^{\alpha_2}\cdots x_{t_n}^{\alpha_n}$, where $t_i\ne t_{i+1}$ for every $i$.
 If the indices are in ascending order, then $w\in R$. Otherwise, we may use the identity $xy=yx+\com x,y$ to
 write $w$ as the sum of a monomial with fewer indices out of order and an element of $\set \com x_1,{x_2}\endset ^T$.
 It follows now by induction that $w\in R+S+T^{(3)}$. 
\end{proof}

\begin{definition}
 A product term $e_{i_1}e_{i_2}\cdots e_{i_n}$ in $\nonunitgrass$ is said to be {\it even} if $n$ is even, otherwise the
 product term is said to be {\it odd}. $u\in \nonunitgrass$ is said to be {\it even} if $u$ is a linear combination of even product
 terms, while $u$ is said to be {\it odd} if $u$ is a linear combination of odd product terms. Let $C$ denote
 the set of all even elements of $\nonunitgrass$, and let $H$ denote the set of all odd elements of $\nonunitgrass$.
\end{definition}

 Note that $C$ and $H$ are subspaces of $\nonunitgrass$, with $C$ closed under multiplication,
 $H^2\subseteq C$, and $CH=HC\subseteq H$. Evidently, $\nonunitgrass=C\oplus H$ as vector spaces.

\begin{lemma}\label{lemma: useful}
  \begin{list}{(\roman{parts})}{\usecounter{parts}}
  \item 
    $C=C_{\nonunitgrass}$. 
  \item 
    For $h,u\in H$, $hu=-uh$. In particular, $h^2=0$ (since $p\ne2$).
  \item 
    Let $g\in \nonunitgrass$, so there exist (unique) $c\in C$ and $h\in H$ such that $g=c+h$. For any 
  positive integer $n$, $g^n=c^n+nc^{n-1}h$. 
  \item 
   If $p>2$, then for $g\in \nonunitgrass$, $g^p=0$.
  \item 
   Let $c_1,c_2\in C$ and $h_1,h_2\in H$, and set $g_1=c_1+h_1$, $g_2=c_2+h_2$. Then
   for any nonnegative integers $m_1,m_2$, $\com g_1,{g_2}g_1^{m_1}g_2^{m_2}=2c_1^{m_1}c_2^{m_2}h_1h_2$
   (where $g_i^0$ and $c_i^0$ are understood to mean that the factors $g_i^0$ and $c_i^0$ are omitted).
  \item 
   Let $u\in \nonunitgrass$. Then $u^{n+1}=0$, where $n$ is the number of distinct basic product terms in the expression
   for $u$ as a linear combination of elements of $\mathcal{B}$.
 \end{list}
\end{lemma}

\begin{proof}
 That $C\subseteq C_{\nonunitgrass}$ and $hu=-uh$ for all $h\in H$, $u\in\kzerox$, it suffices to observe that 
 $e_k(e_{i_1}e_{i_2}\cdots e_{i_n})=(-1)^ne_{i_1}e_{i_2}\cdots e_{i_n}e_k$
 for any $e_k, e_{i_1},e_{i_2},\ldots, e_{i_{n}}$. It remains to prove that $C_{\nonunitgrass}\subseteq C$. Let $u\in C_{\nonunitgrass}$. Then
 $u=c+h$ for some $c\in C$ and $h\in H$. Then $h=u-c\in C_{\nonunitgrass}$ as well. But $h=\sum_{i=1}^n \alpha_iu_i$,
 where for each $i$, $\alpha_i\in k$ and $u_i\in H\cap \mathcal{B}$. Let $j$ be such that $e_j$ does not appear in
 any $u_i$. Since $h\in H$, $he_j=-e_jh$, but since $h\in C_{\nonunitgrass}$, $he_j=e_jh$. Thus $0=2he_j=\sum_{i=1}^n 2\alpha_iu_ie_j$.
 Since $e_j$ does not appear in any $u_i$, we conclude that $2\alpha_i=0$ for every $i$. Thus $2h=0$, and since $p\ne2$,
 it follows that $h=0$ and so $u\in C$. 
 
 For (iii), note that since $c\in C$ is central, $(c+h)^n=\sum_{i=0}^n \choice n,i c^ih^{n-i}$. Since $h\in H$, we
 have $h^2=0$, so only the terms with $i=n-1$ and $n$ survive, so $(c+h)^n=\choice n,{n-1} c^{n-1}h + \choice n,n c^n=
 c^n+nc^{n-1}h$, as required.
 
 For (iv), suppose that $p>2$ and write $c\in C$ as a linear combination of even product terms, say $c=\sum_{j=1}^k \alpha_j u_j$,
 where each $u_j$ is an even product term, hence central. Then 
 \[
  c^p=\sum_{\substack{0\le i_1,i_2,\ldots,i_k\le k\\i_1+i_2+\cdots +i_k=p}}
 \choice p,{i_1,i_2,\ldots,i_k} (\alpha_1u_1)^{i_1}\cdots (\alpha_ku_k)^{i_k}.
 \]
 Since $p$ divides $\choice p,{i_1,i_2,\ldots,i_k}$
 if every $i_j\ne p$, it follows that $c^p=\sum_{i=1}^k \alpha_i^p u_i^p$. Finally, since each $u_i$ is a product term and thus $u_i^2=0$
 for each $i$, it follows that $c^p=0$. For $g\in \nonunitgrass$, we have $g=c+h$ for some $c\in C$ and $h\in H$, and by
 (iii), $g^p=c^p$, so $g^p=0$.
 
 For (v), we have 
 $$
  \com g_1,{g_2}g_1^{m_1}g_2^{m_2}=\com {c_1+h_1},{c_2+h_2}c_1^{m_1}c_2^{m_2}
  =\com h_1,{h_2}c_1^{m_1}c_2^{m_2}=2h_1h_2c_1^{m_1}c_2^{m_2}.
 $$ 
 
 Finally, let $u\in \nonunitgrass$ with $u=\sum_{i=1}^n \alpha_iu_i$, where each $u_i$ is a basic product term. Then 
 $u^{n+1}$ is a linear combination of product terms of the form $u_{i_1}u_{i_2}\cdots u_{i_{n+1}}$, where
 $i_1,i_2,\ldots,i_{n+1}\in \set 1,2,\ldots,n\endset$, so in each product term, there is at least one index $i$
 such that $u_i$ appears twice in the product term, and so each product term is 0. 
\end{proof}  
  
\begin{definition}\label{definition: support}
 For $u=e_{i_1}e_{i_2}\cdots e_{i_n}\in \mathcal{B}$, let $s(u)=\set e_{i_1},e_{i_2},\ldots,e_{i_n}\endset$. Then for any
 $g\in \nonunitgrass$, if $g\ne 0$, let $s(g)=\bigcup_{i=1}^m s(g_i)$, where $g=\sum_{i=1}^m a_ig_i$ with $a_i\in k^*$ and $g_i\in \mathcal{B}$,
 while $s(0)=\nullset$. We shall refer to $s(u)$ as the {\em support} of $u$.
\end{definition}

The next lemma is an adaptation of Lemma 2.9 in \cite{Si}. 
 
\begin{lemma}\label{lemma: nice stuff}
 Let $n$ and $N$ be positive integers and let 
 $w=\sum_{j=N+1}^{N+n} e_{2j-1}e_{2j}\in C$, and let $v=w+e_{2N+2n+1}$. Then the following hold.
  \begin{list}{(\roman{parts})}{\usecounter{parts}}
  \item
    Let $t=e_{2N+1}e_{2N+2}\cdots e_{2N+2n}$. Then $t\in C\cap \mathcal{B}$ and $w^n=n!t\in \gen{t}$,\penalty-10000 so $s(w^n)=s(t)=s(w)$,
   while for any $m<n$, 
   \[
    w^m\in\gen{\set u\in \mathcal{B}\mid \text{ $s(u)\subsetneq s(w)$}\endset}.
    \]
  \item 
   Let $t=e_{2N+1}e_{2N+2}\cdots e_{2N+2n}e_{2N+2n+1}$. Then $t\in H\cap \mathcal{B}$ and $v^{n+1}=(n+1)!t\in \gen{t}$,
   so $s(v^{n+1})=s(t)=s(v)$, while for any $m<n+1$, $v^m\in\gen{\set u\in \mathcal{B}\mid \text{ $s(u)\subsetneq s(v)$}\endset}$.
  \end{list} 
\end{lemma}

\begin{proof}
 Since $e_ie_j$ is central in $\nonunitgrass$ for any $i$ and $j$, it follows that for any positive integer $m$,
 \begin{align*}
   w^m&=\sum_{\substack{ 0\le j_1,\ldots,j_n\le m\\ j_1+\cdots +j_n=m}} \choice m,{j_1,j_2,\ldots,j_n}
  \prod_{k=1}^n (e_{i_{2k-1}}e_{i_{2k}})^{j_k}\\
  &=\sum_{\substack{0\le j_1,\ldots,j_n\le 1\\ j_1+\cdots +j_n=m}} \choice m,{j_1,j_2,\ldots,j_n}
  \prod_{k=1}^n (e_{i_{2k-1}}e_{i_{2k}})^{j_k} \quad\text{since $e_i^2=0$}\\
  &=\mkern 20mu m!\mkern -20mu\sum_{\substack{0\le j_1,\ldots,j_n\le 1\\ j_1+\cdots +j_n=m}} 
  \prod_{k=1}^n (e_{i_{2k-1}}e_{i_{2k}})^{j_k}.
 \end{align*}
 The result follows now by Lemma \ref{lemma: useful} since $w\in C$ and $e_{2N+2n+1}\in H$. Note that by Lemma \ref{lemma: useful} (iii)
 and (vi), $v^{n+1}=w^{n+1}+(n+1)w^ne_{i_{2n+1}}=(n+1)w^ne_{i_{2n+1}}$.
\end{proof}  

\begin{definition}\label{definition: def of M}
 Let $t$ and $n$ be positive integers with $t\le n$, and define $M_{t,n}$ to be the set of all elements of the form 
 \begin{gather}
  x_t^{\alpha_t} \prod_{r=1}^s \com x_{j_{2r-1}},{x_{j_{2r}}}x_{j_{2r-1}}^{\beta_{2r-1}}x_{j_{2r}}^{\beta_{2r}} \label{Mdef i}\\ 
 \intertext{or}
  \biggl( \prod_{r=1}^l x_{i_r}^{\alpha_r}\biggr)x_t^{\alpha_t}\prod_{r=1}^s \com x_{j_{2r-1}},{x_{j_{2r}}}x_{j_{2r-1}}^{\beta_{2r-1}}x_{j_{2r}}^{\beta_{2r}}, \label{Mdef ii} \\
 \intertext{or, in the case $t=n$ only, elements of the form}
 \prod_{r=1}^n x_{i_r}^{\alpha_r} \label{Mdef iii} 
 \end{gather}
 as well, subject to the following requirements.
 \begin{list}{(\alph{parts})}{\usecounter{parts}} 
 \item\label{M (a)} 
  $s,l\ge1$.
 \item 
  $1 \le i_1<\cdots<i_{l}<t$ for elements of type (2.2) or (2.3).
  \item 
   $j_1<\cdots<j_{2s}$.
  \item 
   $\{i_1,\ldots,i_l,t\} \cap \{j_1,\ldots,j_{2s}\}=\nullset$.
  \item 
   $\{i_1,\ldots,i_l,t\} \cup\{j_1,\ldots,j_{2s}\}=\{1,\ldots,n\}$.
  \item 
    For every $k$ with $1 \le k \le l$, $1 \le \alpha_m,\alpha_k$, and, if $p>2$, $\alpha_m,\alpha_k\le p-1$.
  \item 
    For every $r$ with $1 \le r \le 2s$, $0\le\beta_r$, and, if $p>2$, $\beta_r\le p-1$.
 \end{list}
\end{definition}

For any positive integers $m$ and $n$ with $m\le n$, $M_{m,n}\subseteq \siderovset$, and so Definition \ref{definition: total order} defines a total order
on $M_{m,n}$.

\begin{lemma}\label{lemma: m fact}
  Let $m,n$ be positive integers with $m\le n$. Then for any $u \in M_{m,n}$, there exist
  $g_1,g_2,\ldots,g_n$ in $\finitenonunitgrass{z}$, the subalgebra of $\nonunitgrass$ that is generated by 
  $\set e_1,e_2,\ldots,e_{z}\endset$, where $z=2(\deg(u)-\lend{u})-1$, such that the following hold:
   \begin{list}{(\roman{parts})}{\usecounter{parts}}
   \item $0\ne u(g_1,g_2,\ldots,g_n)\in \gen{\prod_{r=1}^z e_i}$;
   \item for any $v\in M_{m,n}$ with $u>v$, $v(g_1,g_2,\ldots,g_n)=0$;
   \item if $p>2$, then for any $g\in \nonunitgrass$, $\com g_m,g g_m^{p-1}=0$.
  \end{list}
\end{lemma}

\begin{proof}
 Let $u\in M_{m,n}$. Then
 $u$ has form as shown in either (\ref{Mdef i}) or (\ref{Mdef ii}), or, in the case $m=n$, of the form (\ref{Mdef iii}),
 and we shall give the proof for $u$ of the form (\ref{Mdef ii}), as it will be evident
 in the argument that the forms (\ref{Mdef i}) and (\ref{Mdef iii}) can be handled in a similar way.
 Under this assumption, there exist indices $i_1,\ldots,i_t$, $j_1,\ldots,j_{2s}$ and integers $\alpha_1,\ldots,\alpha_t$, $\alpha_m$,
 and $\beta_1,\ldots,\beta_{2s}$ satisfying the conditions of Definition \ref{definition: def of M} (a)--(f) such that
 \[
  u(x_1,x_2,\ldots,x_n)=\biggl( \prod_{k=1}^t x_{i_k}^{\alpha_k}\biggr)x_m^{\alpha_m}\prod_{k=1}^s
  \com x_{j_{2k-1}},{x_{j_{2k}}}x_{j_{2k-1}}^{\beta_{2k-1}}x_{j_{2k}}^{\beta_{2k}}.
 \] 

 By repeated applications of Lemma \ref{lemma: nice stuff} (i) and (ii), making appropriate choices for the value of $N$ in
 each case (see (i) below), we find that for the values:
 \newenvironment{myitemize}{%
   \renewcommand{\labelitemi}{--}\begin{itemize}}{\end{itemize}}
 \begin{myitemize}  
  \item $g_m=e_{2N_m+2l-1}+\sum_{i=1}^{l}e_{2N_m+(2i-1)}e_{2N_m+2i}$, where $l=\alpha_m-1$;
  \item for each $k$ with $1\le k\le t$, $g_{i_k}=\sum_{i=1}^{\alpha_k}e_{2N_k+(2i-1)}e_{2N_k+2i}$;
  \item for each $k$ with $1\le k\le 2s$, $g_{j_k}=e_{2N_k+2\beta_k+1}+\sum_{i=1}^{\beta_k}e_{2N_k+(2i-1)}e_{2N_k+2i}$;
 \end{myitemize}         
 the following hold:
 \begin{list}{(\roman{parts})}{\usecounter{parts}}
   \item 
     for any $i,j$ with $i\ne j$, $s(g_i)\cap s(g_j)=\nullset$.
   \item 
     there exists $b_m\in H\cap\mathcal{B}$ such that $s(b_m)=s(g_m)$ and
     $0\ne g_m^{\alpha_m}\in \gen{b_m}$, but for any $i>\alpha_m=\deg_{x_m}(u)$, $g_m^i=0$;
   \item 
     for each $k$ with $1\le k\le t$, there exists $b_{i_k}\in
     C\cap \mathcal{B}$ such that $s(b_{i_k})=s(g_{i_k})$ and
     $0\ne g_{i_k}^{\alpha_k}\in \gen{b_{i_k}}$, but for any
     $j>\alpha_k=\deg_{x_k}(u)$, $g_{i_k}^j=0$;
   \item 
     for each $k$ with $1\le k\le 2s$, if
     $r>\beta_k+1=\deg_{x_k}(u)$, then $g_{j_k}^{r}=0$, while
     if $r\ge \deg_{x_k}(u)$, then for any $g\in G$, $\com
     g_k,g g_k^{r}=0$. Moreover, there exists $b_{j_k}\in H\cap
     \mathcal{B}$ such that  $s(b_{j_k})=s(g_{j_k})$, and for
     each $r$ with $1\le r\le s$, $b_{j_{2r-1}}b_{j_{2r}}\in
     C\cap\mathcal{B}$ and $0\ne \com
     g_{j_{2r-1}},{g_{j_{2r}}}g_{j_{2r-1}}^{\beta_{2r-1}}g_{j_{2r}}^{\beta_{2r}}\in\gen{b_{j_{2k-1}}b_{j_{2k}}}$.
 \end{list}     
 
  Observe that since $|s(g_m)|=2\alpha_m-1=2\deg_{x_m}(u)-1$, and for
  each $r=1,2,\ldots,t$, $|s(g_{i_r})|=2\alpha_r=2\deg_{x_{i_r}}(u)$,
  while for each $r=1,2,\ldots,2s$,
  $|s(g_{j_r})|=2\beta_r+1=2\deg_{x_{j_r}}(u)-1$, we may choose the
  various $N_i's$ so that  $s(g_1)\cup \cdots s(g_n)=\set
  e_1,e_2,\ldots,e_z\endset$, where $z=2\deg_{x_m}(u)-1+\sum_{r=1}^t
  2\deg_{x_{i_r}}(u)+\sum_{r=1}^{2s}(2\deg_{x_{j_r}}(u)-1)=2\deg(u)-2s-1$.

  The first assertion of the lemma now follows immediately, since
  \[
    0\ne u(g_1,g_2,\ldots,g_n)\in\gen{\bigl(\prod_{k=1}^t b_{i_k}\bigr)b_m
                \bigl(\prod_{k=1}^s b_{j_{2k-1}}b_{j_{2k}}\bigr)}
  \] 
  and $\bigl(\prod_{k=1}^t b_{i_k}\bigr)b_m\bigl(\prod_{k=1}^s
  b_{j_{2k-1}}b_{j_{2k}}\bigr)=\prod_{r=1}^z e_r$.

  For the second assertion, suppose that $v\in M_{m,n}$ with $u>v$. Note
  first that if $u>v$ is due to either condition (ii) or condition (iv)
  of Definition \ref{definition: total order}, then there exists an
  index $i$ such that $x_i$ is in the beginning of $u$ and the end of
  $v$ (see the remark following Definition \ref{definition: total
  order}). In such a case, $i\ne m$ since $x_m$ appears in the beginning
  of every element of $M_{m,n}$, and so $i\in \set
  i_1,i_2,\ldots,i_t\endset$. But then $g_i\in C$ and thus
  $v(g_1,g_2,\ldots,g_n)=0$. 

  Next, we observe that if $u>v$ due to condition (i) or condition (iii)
  of Definition \ref{definition: total order}, then there exists $i$
  such that $\deg_{x_i}(v) > \deg_{x_i}(u)$. But then either $i=i_m$ or
  $i=i_k$ for some $k$ with $1\le k\le t$, or else there exist $k$ with $i\in\set
  j_{2k-1},j_{2k}\endset$. If $x_i$ appears in the beginning of $v$,
  then $g_i^{\deg_{x_i}(v)}$ is a factor of $v(g_1,g_2,\ldots,g_n)$, and
  by (ii) or (iii) above,  $g_i^{\deg_{x_i}}(v)=0$ and so
  $v(g_1,g_2,\ldots,g_n)=0$. On the other hand, suppose that $x_i$
  appears in the end of $v$. Then $i\ne m$ and there exist $k,l$ such
  that $i\in \set k,l\endset$ and  $\com
  g_{k},{g_{l}}g_{k}^{(\deg_{x_{k}} u)-1}g_{l}^{\deg_{x_{l}}(v)-1}$ is a
  factor of $v(g_1,g_2,\ldots,g_n)$. If $x_i$ appears in the beginning
  of $u$, then $g_i$ is central and $v(g_1,g_2,\ldots,g_n)=0$.
  Otherwise, $x_i$ appears in the end of $u$, and then by (iii) above,
  we know  that $\com g_{k},{g_{l}}g_{k}^{(\deg_{x_{k}}
  u)-1}g_{l}^{\deg_{x_{l}}(u)-1}=0$ and hence $v(g_1,g_2,\ldots,g_n)=0$,
  as required.
  
  For the third assertion, assume that $p>2$. We have $g_m=c_m+h_m$, where
  $c_m=\sum_{i=N+1}^{N+\alpha_m-1}e_{2i-1}e_{2i}\in C$ and
  $h_m=e_{2N+2\alpha_m-1}\in H$, where $\alpha_m\le p-1$. Thus $c_m$ has
  $\alpha_m-1\le p-2$ (even, hence commuting) summands of the form
  $e_{2i-1}e_{2i}$, and so $c_m^{p-1}=0$. By Lemma \ref{lemma: useful}
  (v), for $g\in \nonunitgrass$, with $g=c+h$, $c\in C$, $h\in H$, we
  have  $\com g_m,g g_m^{p-1}=2c_m^{p-1}h_mh=0$, as required.  
\end{proof}

\begin{definition}
 A polynomial $f\in \konex$ is said to be {\it  essential in its variables} if every variable that appears in any monomial of $f$ appears
 in every monomial of $f$.
\end{definition}

For $f\in \kzerox$, let $M(f)$ denote the number of monomials in $f$.

\begin{lemma}\label{lemma: basic essential}
 Let $V$ be a $T$-space in $\kzerox$.
 \begin{list}{(\roman{parts})}{\usecounter{parts}}
  \item Let $f\in V$. If $f$ is not essential in its variables, then there exist $f_0,f_1$ in $V$ such
        that $f=f_0+f_1$ and $M(f_0)<M(f)$, $M(f_1)<M(f)$.
  \item Let $E_V=\set f\in V\mid f\text{ is essential in its variables}\endset$. Then $V=\gen{E_V}$; that 
        is, $V$ is the linear span of $E_V$.        
  \end{list}
\end{lemma}

\begin{proof}
 We first consider (i).
 Suppose that $f$ is not essential in its variables. Then there exists a variable $x$ that appears in some but not all monomials of $f$. Let
 $f_0$ denote the sum of all monomial terms of $f$ in which $x$ does not appear, and let $f_1=f-f_0$. Since
 $f_0=f\rest{x=0}\in V$, it follows that $f_1\in V$ as well. Furthermore, we have $M(f_0)<M(f)$ and $M(f_1)<M(f)$, as required.
 
  Next, consider (ii).
 If $V=\set 0\endset$, then $E_V=\nullset$ and the result holds. Suppose that $V\ne \set 0\endset$. We prove
 the result by induction on $M(f)$. Of all elements of $V$, let $f$ be one for which $M(f)$ is least possible.
 We note that $M(f)\ge 1$. If $f$ is not essential in its variables, then by Lemma \ref{lemma: basic essential} (i), there exist
 $f_0,f_1\in V$ with $f=f_0+f_1$ and $M(f_0),M(f_1)<M(f)$, which is not possible. Thus $f\in E_V$. Now suppose
 that $f\in V$ and for all $g\in V$ with $M(g)<M(f)$, $g\in \gen{E_V}$. If $f\notin E_V$, then by Lemma \ref{lemma: basic essential} (i),
 there exist $f_0,f_1\in V$ with $M(f_0),M(f_1)<M(f)$ and $f=f_0+f_1$. By our induction hypothesis, $f_0,f_1\in
 \gen{E_V}$, so $f\in \gen{E_V}$ as well. The result follows now by induction.
\end{proof}

\begin{lemma}\label{lemma: identities for inf dim nonunitary}
 If $p=0$, then $T(\nonunitgrass)=T^{(3)}$, while if $p>2$, 
 $T(\nonunitgrass)=\set x_1^p\endset^T+T^{(3)}$.
\end{lemma} 

\begin{proof}
 For $p=0$, this result appears as Corollary 2 of \cite{AK}, while for  $p>2$, it is 
 Theorem 3 of \cite{Si}.
\end{proof}

We are now ready for the main result for the infinite dimensional nonunitary Grassmann algebra over an arbitrary field
of characteristic $p>2$.

\begin{theorem}\label{theorem: main theorem for nonunitary}
  For $k$ any field of characteristic $p\ne2$, $CP(\nonunitgrass)=S+T(\nonunitgrass)$.
\end{theorem}

\begin{proof} Let $U=S+T(\nonunitgrass)$. 
 It is immediate from Lemma \ref{lemma: s central} that $U\subseteq CP(\nonunitgrass)$. 

 For the converse, we note that $CP(\nonunitgrass)$ is a $T$-space, and
 thus by Lemma \ref{lemma: basic essential} (ii), $CP(\nonunitgrass)$ is
 the linear span of its essential polynomials. It suffices therefore to
 prove that any essential element of $CP(\nonunitgrass)$ belongs to $U$.
 So let $f\in CP(\nonunitgrass)$  be essential in $CP(\nonunitgrass)$,
 and suppose that $f\notin U$. Then by Lemma \ref{lemma: rep of
 kzerox},  $f \cong \sum_{j=1}^{t}\alpha_j u_j \mod{U}$, where  for each
 $j$, $\alpha_j\in k^*$, $u_j \in R$, and, if $p>2$, any variable in the beginning of
 $u_j$ has degree at most $p-1$, while any variable in the end of $u_j$ has degree at
 most $p$ (since then $x_1^p\in T(\nonunitgrass)$). Without loss of generality, we
 may assume that the variables that appear in $f$ are
 $x_1,x_2,\ldots,x_n$ for some positive integer $n$. We may further assume that
 for each $j$, $u_j$ is essential in $\kzerox$ in the variables
 $x_1,x_2,\ldots,x_n$. For suppose to the contrary that for some $i$ and
 $j$, $x_i$ does not appear in $u_j$. Let $f_0=\sum_{x_i\text{ not in }u_r}
 \alpha_ru_r$ and $f_1=f-f_0$. Since $x_i$ appears in every monomial
 in $f$ and $f_1$, we have $0\cong f\rest{x_i=0}=f_0+(f_1\rest{x_i=0})
 \cong f_0\mod{U}$, and so $f\cong f_1\mod{U}$.

 Next, we observe that since each $u_j \in R$, $u_j$ has a nonempty beginning. Let  
 \begin{align*}
    m=\max\{\, i\mid 1\le i\le n\ &\text{ and there exists $j$ such}\\
    &\text{that $x_i$ appears in the beginning of $u_j$}\,\}.
 \end{align*}
 Now, for each $j$ such that $x_m$ appears in the beginning of $u_j$,
 $u_j\in M_{m,n}$. Let $f_m$ denote the sum of these terms, so that
 $f_m\in \gen{M_{m,n}}$, and let $f_e=f-f_m$. Suppose that $f_e\in U$.
 Then $f_m\cong f\not\cong 0\mod{U}$, and so $f_m\in M_{m,n}-\set
 0\endset$, which implies by Lemma \ref{lemma: m fact} that $f_m\notin
 CP(\nonunitgrass)$ and thus $f\notin CP(\nonunitgrass)$. Since this is
 not the case, it follows that $f_e\notin U$.  Consider $j$ such that 
 $\alpha_j u_j$ is a summand of $f_e$. Suppose that $\deg_{x_m} u_j<p$,
 so there exist $0<i_1<i_2<\cdots <i_t<m$, positive integers
 $\alpha_1,\ldots,\alpha_t$, $0<j_{1}<\cdots<j_{2s}$, nonnegative
 integers $\beta_1,\ldots,\beta_{2s}$ such that $m\in\set
 j_1,\ldots,j_{2s}\endset$, $\beta_m\le p-2$, and
 \[
   u_j=\prod_{k=1}^t x_{i_k}^{\alpha_k}\,\prod_{r=1}^s\com x_{j_{2r-1}},{x_{j_{2r}}} 
                             x_{j_{2r-1}}^{\beta_{2r-1}}x_{j_{2r}}^{\beta_{2r}}.
 \]

 We may assume without loss of generality that $m$ is odd, since the
 argument for the case when $m$ is even can be converted,  modulo $T^{(3)}$,
 to the case where $m$ is odd by a sign change in $\alpha_j$. Thus
 $m=2d-1$ for some $d$ with $1\le d\le s$, and for convenience, let
 $\beta=\beta_{j_{2d-1}}$. We apply Corollary \ref{corollary: need for R
 gen}, taking $v=\prod_{k=1}^t x_{i_k}^{\alpha_k}$ and
 $u=\prod_{\substack{r=1\\r\ne d}}^s
             \com x_{j_{2r-1}},{x_{j_{2r}}} x_{j_{2r-1}}^{\beta_{2r-1}}x_{j_{2r}}^{\beta_{2r}}$, to find that
 \begin{align*}
  u_j&\overset{S+T^{(3)}}{\cong} (\beta+1)^{-1}\com x_{j_{2d}},v x_m^{\beta+1}x_{j_{2d}}^{\beta_{2d}}u\\
    &\overset{T^{(3)}}{\cong}(\beta+1)^{-1}\sum_{k=1}^t\alpha_k\com x_{j_{2d}},{x_{i_k}} 
      x_{i_k}^{\alpha_k-1}\bigl(\prod_{\substack{l=1\\l\ne k}}^tx_{i_l}^{\alpha_l}\bigr)
      x_m^{\beta+1}x_{j_{2d}}^{\beta_{2d}}u\quad\text{\begin{tabular}[t]{l} by Lemma \ref{lemma: handy}\\ (iii) and (vi)\end{tabular}}\\
    &\overset{T^{(3)}}{\cong}\sum_{k=1}^t(\beta+1)^{-1}\alpha_k
      \bigl(\prod_{\substack{l=1\\l\ne k}}^tx_{i_l}^{\alpha_l}\bigr)x_m^{\beta+1}
      \com x_{j_{2d}},{x_{i_k}} x_{j_{2d}}^{\beta_{2d}}x_{i_k}^{\alpha_k-1}u
      \quad\text{by Lemma \ref{lemma: handy} (vi)}
 \end{align*}
 Now for each $k$, $\displaystyle\bigl(\prod_{\substack{l=1\\l\ne k}}^t
 x_{i_l}^{\alpha_l}\bigr)x_m^{\beta+1}\com x_{j_{2d}},{x_{i_k}}
 x_{j_{2d}}^{\beta_{2d}}x_{i_k}^{\alpha_k-1}u$ is equal to 
 \[
 \bigl(\prod_{\substack{l=1\\l\ne k}}^tx_{i_l}^{\alpha_l}\bigr)x_m^{\beta+1}\com x_{j_{2d}},{x_{i_k}} 
  x_{j_{2d}}^{\beta_{2d}}x_{i_k}^{\alpha_k-1} \prod_{\substack{r=1\\r\ne d}}^s
  \com x_{j_{2r-1}},{x_{j_{2r}}} x_{j_{2r-1}}^{\beta_{2r-1}}x_{j_{2r}}^{\beta_{2r}},
 \]
 and by Lemma \ref{lemma: handy}, working modulo $T^{(3)}$, the end of this
 element can be rearranged so as to give an element of $R$ with
 beginning $\bigl(\prod_{\substack{l=1\\l\ne k}}^tx_{i_l}^{\alpha_l}\bigr)x_m^{\beta+1}$.	       
 We have proven now that if $p=0$ or else $p>2$ and $\deg_{x_m} u_j<p$, then $u_j$ is congruent to a linear
 combination of elements of $M_{m,n}$. It follows that $m$ is such that
 there exist $f_m\in \gen{M_{m,n}}$ and, if $p>2$,
 \begin{align*}
   f_e\in \left<\{\, u\in R\mid \deg_{x_m}u=\right.&p,\ \text{ and for every $i$, if $x_i$}\\
    &\left.\text{appears in the beginning of $u$, then $i<m$}\,\}\right>  
 \end{align*}
 while if $p=0$, take $f_e=0$,
 such that $f\cong f_m+f_e\mod{S+T^{(3)}}$. Suppose now that $m$ is minimal
 with respect to this property. If $f_m\cong 0\mod{S+T^{(3)}}$, then $f\cong
 f_e\mod{S+T^{(3)}}$, which contradicts our choice of $m$ since $f_e$ is a
 linear combination of elements of $R$ in whose beginning only elements
 $x_i$ with $i<m$ appear. Thus $f_m\in \gen{M_{m,n}}-\set0\endset$, and
 so by Lemma \ref{lemma: m fact},  there exist $g_1,g_2,\ldots,g_n\in
 \nonunitgrass$ such that $0\ne f_m(g_1,g_2,\ldots,g_n)\in\gen{e_1e_2\cdots e_z}$
 where $z=2(\deg(u)-\lend{u})-1$ and for any $g\in \nonunitgrass$, 
 $\com g_m,g g_m^{p-1}=0$. Since $z$ is odd, $f_m(g_1,g_2,\ldots,g_n)\notin
 C_{\nonunitgrass}$. If $p>2$, then $f_e=\sum
 \gamma_iv_i$, where for each $i$, $v_i\in R$ and $\deg_{x_m}v_i=p$, from which 
 we obtain that in each $v_i$, $x_m$ appears in a term of the form $\com
 x_m,yx_m^{p-1}y^{\beta}$, and so $v_i(g_1,g_2,\ldots,g_n)=0$ for each
 $i$. Thus if $p>2$, $f_e(g_1,g_2,\ldots,g_n)=0$, while if $p=0$, this holds
 by definition. We therefore have $f(g_1,g_2,\ldots,g_n)\cong
 f_m(g_1,g_2,\ldots,g_n)+0\mod{C_{\nonunitgrass}}$. Since 
 $f_m(g_1,g_2,\ldots,g_n)\notin C_{\nonunitgrass}$, this implies that $f\notin
 CP(\nonunitgrass)$, contrary to fact. Since this contradiction follows from our assumption
 that $f\notin U$, it follows that $f\in U$,  as required.
\end{proof}

\begin{corollary}\label{t space generators for cp nonunitary}
 If $k$ is a field of characteristic zero, then 
 $$
  CP(\nonunitgrass)=\set \com x_1,{x_2},\com x_1,{x_2}\com x_3,{x_4}\endset^S,
 $$
 while if $k$ has characteristic $p>2$, then
 $$
  CP(\nonunitgrass)=\set \com x_1,{x_2},x_1^p,x_2x_1^p,\com x_1,{x_2}\com x_3,{x_4}\endset^S+\set w_n\mid  n\ge 1\endset^S,
 $$ 
 where for each $n\ge1$, $w_n=\prod_{k=1}^n\com {x_{2k-1}},{x_{2k}}x_{2k-1}^{p-1}x_{2k}^{p-1}$.
\end{corollary} 


\section{The central polynomials of the infinite dimensional unitary Grassmann
           algebra over an infinite field of characteristic $p>2$}  


 Note that $\unitgrass=\gen{1}\oplus \nonunitgrass=\gen{1}\oplus C\oplus H$, and
 $C_{\unitgrass}=\gen{1}\oplus C$.
 For each  $g\in \unitgrass$, $g=\alpha+c+h$ for some $\alpha\in k$,
 $c\in C$, and $h\in H$. 

\begin{lemma}[\cite{Gi}, Theorem 6]\label{lemma: giambruno's identities}
 $T(\unitgrass)=T^{(3)}$.
\end{lemma}

\begin{definition}\label{definition: def of S1}
 If $p=0$, let $S_1=S$, while if $p>2$, let 
 \[
   S_1=\set \com x_1,{x_2},x_1^p\endset^S+\set x_{1}^p
    \prod_{i=1}^k\com x_{2i},{x_{2i+1}}x_{2i}^{p-1}x_{2i+1}^{p-1}\mid k\ge1\endset^S.
 \]    
\end{definition}

Note that if $p>2$, $S\subseteq S_1$, by virtue of the substitution $x_1\to 1$, while if
$p=0$, then $S=S_1$. Thus in every case, $S\subseteq S_1$.

\begin{lemma}\label{lemma: for unitary}
 Let $p>2$ and $u\in \siderovset$ with $\lbeg{u}>0$. 
 If $\deg_{x_i}(u)\cong0\mod{p}$ for every $x_i$ that appears in the beginning
 of $u$, then $u\in S_1+T^{(3)}$.
\end{lemma}

\begin{proof} By Lemma~\ref{lemma: handy}~(vii),
 $x_{i_1}^px_{i_2}^p\cdots x_{i_t}^p\cong (x_{i_1}\cdots
 x_{i_t})^p\mod{T^{(3)}}$, and so it suffices to prove the result for the case
 $$
   u=\prod_{k=1}^t x_{i_k}^{\alpha_k} \prod_{k=1}^s 
 \com x_{j_{2i-1}},{x_{j_{2i}}}x_{j_{2i-1}}^{\beta_{2i-1}}{x_{j_{2i}}}^{\beta_{2i}}
 $$
 with $t=1$. Thus we consider $u= x_{i_1}^{\alpha_1} \prod_{k=1}^s \com
 x_{j_{2i-1}},{x_{j_{2i}}}x_{j_{2i-1}}^{\beta_{2i-1}}{x_{j_{2i}}}^{\beta_{2i}}\in\siderovset$
 with $\alpha_1\cong0\mod{p}$. Since $x^p$
 is central in $\konex$ modulo $T^{(3)}$, it follows from Lemma \ref{lemma:
 handy} (vii) that we may assume  that $\beta_i\le p-1$ for each $i$. If
 $\beta_i=p-1$ for each $i$, then $u\in S_1$, so we need only consider
 the situation when there exists $i$ with $\beta_i<p-1$. Without loss of
 generality, we may assume that $i$ is odd, say $i=2r-1$ for some $r$
 with $1\le r\le s$, so we have $1\le 1+\beta_{2r-1}\le p-1$. Let $v=
 \prod_{\substack{k=1\\k\ne r}}^s \com x_{j_{2i-1}},{x_{j_{2i}}}
 x_{j_{2i-1}}^{\beta_{2i-1}}{x_{j_{2i}}}^{\beta_{2i}}$, so by Lemma 
 \ref{lemma: handy} (vi), $v$ is central modulo $T^{(3)}$. Using Lemma 
 \ref{lemma: handy} (i) and (iii), working modulo $T^{(3)}$, we obtain that 
 \begin{align*}
  \com x_{i_1}^px_{j_{2r-1}}^{1+\beta_{2r-1}}{x_{j_{2r}}^{\beta_{2r}}v},{x_{j_{2r}}}&\\
   &\hskip-70pt\cong x_{i_1}^p\com x_{j_{2r-1}}^{1+\beta_{2r-1}}{x_{j_{2r}}^{\beta_{2r}}v},{x_{j_{2r}}}
   +p\com x_{i_1},{x_{j_{2r}}}x_{i_1}^{p-1}x_{j_{2r-1}}^{1+\beta_{2r-1}}{x_{j_{2r}}^{\beta_{2r}}v}\\
   &\hskip-70pt\cong(1+\beta_{2r-1})x_{i_1}^p\com x_{j_{2r-1}},{x_{j_{2r}}}x_{j_{2r-1}}^{\beta_{2r-1}} x_{j_{2r}}^{\beta_{2r}}v +
   x_{i_1}^p\com {x_{j_{2r}}^{\beta_{2r}}v},{x_{j_{2r}}}x_{j_{2r-1}}^{1+\beta_{2r-1}} \\
   &\hskip-70pt\cong(1+\beta_{2r-1})u.
 \end{align*}
 
 Since $\com x_{i_1}^px_{j_{2r-1}}^{1+\beta_{2r-1}}{x_{j_{2r}}^{\beta_{2r}}v},{x_{j_{2r}}}\in S_1$ 
 and $1+\beta_{2r-1}\not\cong0\mod{p}$, it follows that $u\in S_1+T^{(3)}$.
\end{proof} 
 
\begin{definition}\label{definition: def of R1}
 Let $R_1$ be the subspace of $\konex$ that is spanned by  $1$ and
 \begin{align*}
  \{\, u\in \siderovset \mid \lbeg{u}>0,\ &\text{and, if $p>2$, for some $x_i$ in the}\\
   &\text{beginning of $u$, $\deg_{x_i}(u)\not\cong 0\mod{p}$}\,\}.
  \end{align*}
\end{definition} 
 
\begin{corollary}\label{corollary: rep of k1x}
 $\konex=R_1+S_1+T^{(3)}$.
\end{corollary}

\begin{proof}
 In the case $p=0$, $R_1=k+R$ and $S_1=S$, so by Lemma \ref{lemma: rep of kzerox}, we have
 $\konex=k+\kzerox=k+R+S+T^{(3)}=R_1+S_1+T^{(3)}$. Suppose now that $p>2$. Our proof in this
 case is similar to that of Lemma \ref{lemma: rep of
  kzerox}. By Proposition \ref{proposition: 2-comm ideal} and Corollary \ref{corollary: alt description of S gen},
 each element of $\set \com x_1,{x_2}\endset^T$ 
 is in $R+S_1+T^{(3)}$. Moreover, any $u\in R$ in which every variable that appears in the beginning of $u$ has
 degree congruent to 0 modulo $p$ belongs to $S_1+T^{(3)}$ by Lemma \ref{lemma: for unitary}.
 Thus $\set \com x_1,{x_2}\endset^T\subseteq R_1+S_1+T^{(3)}$. 
 Since $1\in R_1$ (and $S_1$ in this case), it therefore suffices to prove that every monomial belongs to $R_1+S_1+T^{(3)}$.
 From this point on, proceed as in the proof of Lemma \ref{lemma: rep of kzerox}.
\end{proof}
  
\begin{definition}\label{definition: def of Mmnprime}
 For positive integers $m,n,r_1,r_2,\ldots,r_n$ with $m\le n$ and
 $r_m\not\cong 0\mod{p}$, let $M'_{m,n}(r_1,r_2,\ldots,r_n)$ denote the
 set of all $u\in\siderovset$ with $\lbeg{u}>0$ for which $x_m$ appears 
 in the beginning of $u$ and $\deg_{x_i} u = r_i$ for every $i$.
\end{definition}

 It is useful to observe that the conditions of Definition
 \ref{definition: total order} also define a total ordering on
 $M'_{m,n}(r_1,r_2,\ldots,r_n)$.

\begin{lemma}\label{lemma: m' fact}
  Let $m,n,r_1,r_2,\ldots,r_n$ be positive integers with $m\le n$ and, if $p>2$, $r_m\not\cong 0\mod{p}$.
  Then for any $u \in M'_{m,n}(r_1,\ldots,r_n)$, there exist
  $g_1,g_2,\ldots,g_n$ in $\finiteunitgrass{z}$, the subalgebra of $\unitgrass$ that is generated by 
  $\set e_1,e_2,\ldots,e_{z}\endset$, where $z=2\lend{u}+1$, such that the following hold:
   \begin{list}{(\roman{parts})}{\usecounter{parts}}
   \item $u(g_1,g_2,\ldots,g_n)$ has nonzero odd part in $\gen{\prod_{r=1}^z e_i}$;
   \item for any $v\in M'_{m,n}(r_1,\ldots,r_n)$ with $u>v$, $v(g_1,g_2,\ldots,g_n)=0$;
  \end{list}
\end{lemma}

\begin{proof}
  Let $u=\prod_{l=1}^t x_{i_l}^{r_{i_l}} \prod_{l=1}^s \com x_{j_{2l-1}},{x_{j_{2l}}}
           x_{j_{2l-1}}^{r_{j_{2l-1}}-1}x_{j_{2l}}^{r_{j_{2l}}-1}$,
  and suppose that $k$ is such that $i_k=m$. Let $g_{m}=1+e_{m}$, and for $l\ne k$, let
  $g_{i_l}=1$. Let $g_{j_l}=1+e_{j_l}$ for every $l$ with $1\le l\le
  2s$. Then
  \[
    u(g_1,g_2,\ldots,g_n)=(1+r_{m}e_{m})2^s\prod_{l=1}^s e_{j_{2l-1}}e_{j_{2l}}
  \]
  has nonzero (since $r_m\not\cong0\mod{p}$ if $p>2$) odd part
  $r_{m}e_{m}2^s\prod_{l=1}^s e_{j_{2l-1}}e_{j_{2l}}\in \gen{\prod_{r=1}^{2s+1} e_i}$.

  For any $v\in M'_{m,n}(r_1,\ldots,r_n)$ with $u>v$, we have $\lend{v}\ge \lend{u}$ and thus there
  exists an index $l\ne m$ with $x_l$ in the beginning of $u$ and in the
  end of $v$. But then $g_l=1$ and so $v(g_1,g_2,\ldots,g_n)=0$. 
\end{proof}

\begin{corollary}\label{corollary: need for main unitary result}
 For any positive integers $m,n,r_1,r_2,\ldots,r_n$ with $m\le n$ and, if $p>2$, $r_m\not\cong 0\mod{p}$, 
 \[
   \gen{M'_{m,n}(r_1,\ldots,r_n)}\cap CP(\unitgrass)=\set 0\endset.
 \]
\end{corollary}

\begin{proof}
  For convenience, let $M'=M'_{m,n}(r_1,\ldots,r_n)$, and let
  $f\in\gen{M'}$, $f\ne 0$. Then $f=\sum_{i=1}^t \alpha_i u_i$ for 
  $\alpha_i\in k^*$, $u_i\in M'$, with $u_1>u_2>\ldots>u_t$.  By
  Lemma \ref{lemma: m' fact}, there exist $g_1,g_2,\ldots,g_n\in\unitgrass$
  such that $u_1(g_1,g_2,\ldots,g_n)$ has nonzero odd part and for each
  $j>1$, $u_j(g_1,g_2,\ldots,g_n)=0$. Thus
  $f(g_1,g_2,\ldots,g_n)=u_1(g_1,g_2,\ldots,g_n) \notin C_{\unitgrass}$,
  and so $f\notin CP(\unitgrass)$.    
\end{proof}

\begin{definition}\label{definition: multihomogeneous}
 A nonzero polynomial $f\in \konex$ is said to be {\it multihomogeneous
 of type} $(r_1,r_2,\ldots,r_n)$ if $f=\sum_{i=1}^t\alpha_i u_i$, where
 for each $i$, $\alpha_i\ne 0$, $u_i\in k_0\langle x_1,\ldots,x_n\rangle$, 
 and $\deg_{x_j} u_i=r_j$ for each $j$. 
\end{definition}

 It is an immediate consequence of 
 \cite{Ch}, Lemma 1.3 (1), that if $k$ is infinite, then every $T$-space of $\konex$ is 
 generated (as a $T$-space) by its multihomogeneous elements.

\begin{theorem}\label{theorem: central in unitary}
 Let $k$ be an infinite field of characteristice $p>2$. Then
 $CP(\unitgrass)=S_1+T(\unitgrass)$.
\end{theorem}

\begin{proof}
 Let $U_1=S_1+T(\unitgrass)=S_1+T^{(3)}$. By Lemma \ref{lemma: s central}, $S\subseteq CP(\unitgrass)$.
 Suppose first that $p=0$. Then $U_1=S+T^{(3)}$, and so $U_1\subseteq CP(\unitgrass)$. On the other hand, 
 if $p>2$, then since 
 $g^p\in C_{\unitgrass}$ for every $g\in \unitgrass$, it follows that $x_1^p\in CP(\unitgrass)$, which together 
 with the fact that $S\subseteq CP(\unitgrass)$ yields $U_1\subseteq CP(\unitgrass)$. Suppose
 that $CP(\unitgrass)-U_1\ne\nullset$, and let $f\in
 CP(\unitgrass)-U_1$. By the remarks above, if every multihomogeneous
 element of $CP(\unitgrass)$ belonged to $U_1$, then
 $CP(\unitgrass)\subseteq U_1$. We may therefore assume that $f$ is
 multihomogeneous, say of type $(r_1,r_2,\ldots,r_n)$. 
 
 By Corollary \ref{corollary: rep of k1x}, $f\cong\sum_{i=1}^b\alpha_i u_i\mod{U_1}$,
 where $u_i\in R_1$ for each $i$. Let $Y=\set u_i\mid 1\le i\le b,\
 \text{$u_i$ is multihomogeneous of the same type as $f$}\endset$, and
 let $W=\set u_i\mid 1\le i\le b\endset - Y$. Then there exist
 $y_1,y_2\in U_1$, such that $y_1$ is either 0 or else is
 multihomogeneous of the same type as $f$, while $y_2$ is either 0 or
 else is multihomogeneous of type different from that of $f$, and, for 
 $f_1=\sum_{u_i\in Y} \alpha_i u_i$ and $f_2=\sum_{u_i\in W} \alpha_i
 u_i$, then $f=f_1+y_1+f_2+y_2$. Necessarily, $f_1+y_1$ is either 0 or
 else multihomogeneous of the same type as $f$, while $f_2+y_2$ is
 either not multihomogeneous or is multihomogeneous of a type different
 from that of $f$. But then $f-(f_1+y_1)$ is either 0 or is
 multihomogeneous of type that of $f$, so it follows that
 $f-(f_1+y_1)=0$; that is, $f\cong f_1\mod{U_1}$. We may therefore
 assume that $W=\nullset$.

 Observe that for every $j$, $\lbeg{u_j}>0$, and, if $p>2$, there must exist 
 an index $i$ such that $r_i\not\cong 0\mod{p}$ and
 $x_i$ is in the beginning of $u_j$ for some $j$, because
 otherwise $u_j\in U_1$ for all $j$ and so $f\in U_1$, which is not the case. 
 Let 
 \begin{align*}
 m=\max\{\, i \mid &\text{$x_i$ is in the beginning of $u_j$ for some $j$, and,}\\
  &\hskip110pt\text{if $p>2$, $r_i\not\cong 0\mod{p}$}\,\}
 \end{align*}
 
   Suppose that for some  index $j$, 
 $x_m$ appears in the end of $u_j$. Without loss of generality, we may
 assume that
 \[
  u_j=x_{i_1}^{r_{i_1}}x_{i_2}^{r_{i_2}}\cdots
 x_{i_t}^{r_{i_t}}\prod_{k=1}^s\bigl(\com x_{j_{2k-1}},{x_{j_{2k}}}
 x_{j_{2k-1}}^{\beta_{2k-1}}x_{j_{2k}}^{\beta_{2k}}\bigr),
 \]
 where $\beta_i=r_{j_i}-1$ for each $i$,  and $m=j_{2l-1}$ for some $l$
 with $1\le l\le s$. By Corollary \ref{corollary: need for R gen} (since
 $S\subseteq S_1$),
 \begin{align*}
   u_j&\overset{U_1}{\cong} r_m^{-1}\com x_{j_{2l}},{x_{i_1}^{r_{i_1}}x_{i_2}^{r_{i_2}}\cdots x_{i_t}^{r_{i_t}}}x_m^{r_m}x_{j_{2l}}^{\beta_{2l}}
     \prod_{\substack{k=1\\k\ne l}}^s\bigl(\com x_{j_{2k-1}},{x_{j_{2k}}}x_{j_{2k-1}}^{\beta_{2k-1}}x_{j_{2k}}^{\beta_{2k}}\bigr)\\
     &\overset{U_1}{\cong}\sum_{a=1}^t r_m^{-1}r_{i_a}\com x_{j_{2l}},{x_{i_a}}x_{i_a}^{r_{i_a}-1}x_{j_{2l}}^{\beta_{2l}}\bigl(\prod_{\substack{h=1\\h\ne a}}^t x_{i_h}^{r_{i_h}}\bigr)x_m^{r_m}  
     \hskip-2pt\prod_{\substack{k=1\\k\ne l}}^s\bigl(\com x_{j_{2k-1}},{x_{j_{2k}}}x_{j_{2k-1}}^{\beta_{2k-1}}x_{j_{2k}}^{\beta_{2k}}\bigr).\\
     &\overset{T^{(3)}}{\cong}\sum_{a=1}^t r_m^{-1}r_{i_{a}}\bigl(\prod_{\substack{h=1\\h\ne a}}^t x_{i_h}^{r_{i_h}}\bigr)x_m^{r_m}  
     \com x_{j_{2l}},{x_{i_a}}x_{i_a}^{r_{i_a}-1}x_{j_{2l}}^{\beta_{2l}}\hskip-2pt\prod_{\substack{k=1\\k\ne l}}^s\bigl(\com x_{j_{2k-1}},{x_{j_{2k}}}x_{j_{2k-1}}^{\beta_{2k-1}}x_{j_{2k}}^{\beta_{2k}}\bigr).\\
 \end{align*}
 If $p>2$, then since $\com x_k^p,{x_m}\cong p\com x_k,{x_m}x_k^{p-1}=0\mod{T^{(3)}}$, it
 follows that $x_m^{r_m}$ can be moved, modulo $T^{(3)}$, to the left of
 each $x_{i_a}^{r_{i_a}}$ in the beginning of $u_j$ for which $m<i_a$
 (by choice of $m$, $r_{i_a}\cong 0\mod{p}$ if $m<i_a$). Furthermore, by
 application of Lemma \ref{lemma: handy} (iv) and (vi), the end of each
 summand may be manipulated modulo $T^{(3)}$ so as to present $u_j$ as a sum
 of multihomogeneous elements of type $(r_1,r_2,\ldots,r_n)$, each in
 $R_1-U_1$. We may therefore assume that $x_m$ appears in the beginning
 of each $u_j$. But then $\sum_{i=1}^b \alpha_i u_i\in
 \gen{M'_{m,n}(r_1,\ldots,r_n)}$, and since $f\cong \sum_{i=1}^b \alpha_i
 u_i\mod{U_1}$, it follows from Corollary \ref{corollary: need for main unitary
 result} that $\sum_{i=1}^b \alpha_iu_i=0$.  Since this implies that
 $f\in U_1$, which contradicts our choice of $f$,  we must have
 $CP(\unitgrass)-U_1=\nullset$, as required.
\end{proof}

\begin{corollary}
 Let $k$ be an infinite field. If $k$ has characteristic zero, then $CP(\unitgrass)=
 \set \com x_1,{x_2}, \com x_1,{x_2}\com x_3,{x_4}\endset^S$, while if $k$ has characteristic $p>2$, then
  $CP(\unitgrass)=\set \com x_1,{x_2},x_1^p,\com x_1,{x_2}\com x_3,{x_4}\endset^S+\set x_{2n+1}^pw_n\mid  n\ge 1\endset^S$,
  where for each $n\ge1$, $w_n=\prod_{k=1}^n\com {x_{2k-1}},{x_{2k}}x_{2k-1}^{p-1}x_{2k}^{p-1}$.
\end{corollary}


\begin{thebibliography}{00}
\bibitem{AK} A.Z. Anan'in and A.R. Kemer, {\em Varieties of associative algebras whose lattice of subvarieties is distributive},
Siberian Math. J. {\bf 17}, (1976), no. 4, 549--554 (1977).

\bibitem{Br} Brandao, Antonio Pereira, Ph. D. thesis. 

\bibitem{Ch} Chuluundorj Bekh-Ochir and D. M. Riley, {\em On the Grassmann $T$-space}, Journal of Algebra and its Applications {\bf 7} (2008), no. 3, 319--336. 

\bibitem{D} V.S. Drensky, {\em Free algebras and $PI$-algebras}. Graduate course in algebra. Springer-Verlag Singapore, Singapore, 2000.

\bibitem{Fo} E. Formanek, Central polynomials for matrix rings, J. Algebra 23 (1972), 129--132.

\bibitem{Gi} A. Giambruno and P. Koshlukov, {\em On the identities of the Grassmann algebras in characteristic $p>0$}, Isr. J. Math.
{\bf 122} (2001), 305--316.

\bibitem{Ka} I. Kaplansky, Problems in the theory of rings, Report of a Conference on Linear Algebras, June, 1956, National Acad. Sci.-National Research Council,
Washington, Publ. 502 (1957), 1--3.
 

\bibitem{L} V. N. Latyshev, {\em On the choice of a basis in a $T$-ideal}, Siberian  Math.  J. {\bf 4} (1963), no. 5, 1122--1127.

\bibitem{Raz} Yu. P. Razmyslov, On a problem of Kaplansky (Russian), Izv. Akad. Nauk SSSR, Ser. Mat. 37 (1973), 483--501. Translation: Math. USSR, Izv. 7 (1973), 479--496.

\bibitem{Si} Plamen Zh. Chiripov and Plamen N. Siderov, 
{\em On bases for identities of some varieties of associative algebras} 
(Russian), PLISKA Stud. Math. Bulgar, {\bf 2} (1981), 103-115.


\bibitem{At} A. N. Stojanova-Venkova, {\em Bases of identities of Grassmann algebras}, Serdica 6 (1980), no. 1, 63--72.


\end{thebibliography}
\end{document}